\documentclass[a4paper]{article}
\usepackage[english]{babel}
\usepackage[utf8x]{inputenc}
\usepackage[T1]{fontenc}
\usepackage{times}
\usepackage{pgf,tikz,pgfplots}
\usepackage{mathrsfs}
\usepackage{tabularx}
\usepackage[round]{natbib}
\usepackage{amsmath}
\usepackage{algorithm}
\usepackage{algpseudocode}
\usepackage{amsfonts} 
\usepackage{soul}
\usepackage{bm}
\usepackage{amssymb}
\usepackage{multirow}
\usepackage{float}
\usepackage{pdflscape}
\usepackage{afterpage}
\usepackage[colorinlistoftodos]{todonotes}
\usepackage{subcaption}
\usepackage{caption}
\usepackage{booktabs,makecell}
\usepackage{graphicx}
\usepackage{csvsimple}
\usepackage{xcolor}
\usepackage{capt-of}
\usepackage[colorlinks=true, allcolors=blue]{hyperref}
\usepackage{fixmath}
\usepackage{forest}
\pgfplotsset{compat=1.14}
\usepackage{array}
\usepackage{arydshln}

\DeclareMathOperator*{\argmin}{argmin} 

\title{A Two-step Heuristic for the Periodic Demand Estimation Problem}

\author{Greta Laage\textsuperscript{a,b,}\footnote{Corresponding author.}, Emma Frejinger\textsuperscript{c,b}, Gilles Savard\textsuperscript{a,d}}

\date{\today}
\author{Greta Laage\footnote{\'Ecole Polytechnique de Montr\'eal, Montr\'eal, Canada, {greta.laage@polymtl.ca}} \and Emma Frejinger \footnote{Department of Computer Science and Operations Research, Universit\'e de Montr\'eal, Montr\'eal, Canada} \and Gilles Savard\footnote{IVADO and \'Ecole Polytechnique de Montr\'eal,  Montr\'eal, Canada}}

\begin{document}

\maketitle

\begin{abstract}
Freight carriers rely on tactical plans to satisfy demand in a cost-effective way. For computational tractability in real large-scale settings, such plans are typically computed by solving deterministic and cyclic formulations. An important input is the \emph{periodic demand}, i.e., the demand that is expected to repeat in each period of the planning horizon. Motivated by the discrepancy between time series forecasts of demand in \emph{each period} and the periodic demand, \cite{Laage2021} recently introduced the Periodic Demand Estimation (PDE) problem and showed that it has a high value. However, they made strong assumptions on the solution space so that the problem could be solved by enumeration. In this paper we significantly extend their work. We propose a new PDE formulation that relaxes the strong assumptions on the solution space. We solve large instances of this formulation with a two-step heuristic. The first step reduces the dimension of the feasible space by performing clustering of commodities based on instance-specific information about demand and supply interactions. The formulation along with the first step allow to solve the problem in a second step by either metaheuristics or the state-of-the-art black-box optimization solver NOMAD. In an extensive empirical study using real data from the Canadian National Railway Company, we show that our methodology produces high quality solutions and outperforms existing ones.

\end{abstract}


\section{Introduction}

Freight carriers rely on tactical plans to satisfy demand in a cost-effective way. In this context, Service Network Design (SND) is an important class of problems \citep{wieberneit2008service}.
Deterministic and cyclic SND formulations are mostly used for real large-scale applications for the sake of computational tractability and
such formulations rely on an accurate representation of \textit{periodic demand}. That is, the demand expected to repeat in each period (e.g., a week) of the tactical planning horizon (e.g., a few months). A tactical plan minimizing costs while satisfying the periodic demand is  defined over the period and repeated over the planning horizon. 

Time series models are widely used to forecast demand for tactical planning. Motivated by the discrepancy between such time series forecasts of demand in \emph{each period} and the periodic demand, \cite{Laage2021} introduce the Periodic Demand Estimation (PDE) problem. The aim is to map the time series forecasts to a periodic demand minimizing the costs over the tactical planning horizon. Their results show that addressing the problem can lead to substantially reduced costs compared to standard practice, i.e., using forecasts averaged over time. They propose a three-level formulation whose solution is the mapping. Under weak assumptions, the multilevel problem can be solved sequentially to optimality. The computational cost is nevertheless high since the lower levels correspond to large-scale combinatorial non-convex and non-differentiable problems. For the sake of tractability, \cite{Laage2021} restrict the solution space to only a few mappings and solve the PDE problem by enumeration. This strong assumption is the motivation behind our work.

Our objective in this paper is to find high-quality solutions to large-scale instances of the PDE problem. For this purpose we relax the assumptions of \cite{Laage2021} on the PDE problem's solution space. In turn, this leads to a problem that cannot be solved by enumeration. We therefore aim to devise a solution approach that can solve large instances. Moreover, the solution approach should be easy to implement and to use in practice.

This paper offers both methodological and empirical contributions. Firstly, we propose a new formulation based on a continuous extension of the discrete solution space in \cite{Laage2021}. We formulate the periodic demand as a deviation from the average forecasts, and we introduce new variables, referred to as deviation coefficients. 
They allow for an intuitive interpretation of the PDE solution. Secondly, we propose a two-step heuristic. The first step consists in data-driven clustering based on instance-specific demand and supply interactions. It reduces the number of variables through heuristic variable fixing. More precisely, commodities grouped into a same cluster are assumed to have equal deviation coefficients. Our main methodological contributions reside in the new formulation and in the first step of the heuristic. Indeed, these, in turn, allow to leverage state-of-the-art metaheuristics in the second step, where we 
use an off-the-shelf black-box optimization solver \citep[NOMAD][]{le2011algorithm}. We compare its performance to two tailored local search heuristics that we propose as baselines. 
Thirdly, we report results in an extensive empirical study of a real large-scale application from the Canadian National Railway Company, our industrial partner. We compare our results to those in \cite{Laage2021} and show that our solutions are of higher quality than existing ones. Moreover, we show that the clustering step is important to solve large-scale instances.

The remainder of the paper is structured as follows. Next we present related work and 
in Section~\ref{section:methodology}, we introduce our formulation and the two-step heuristic proposed to solve the PDE problem. Then, we outline our large-scale application in Section~\ref{section:application}. 
Finally, we report empirical results in Section~\ref{section:results} and conclude with directions for future research in Section~\ref{section:conclusion}.

\section{Related Work}
\label{section:related_works}
Since this work builds on \cite{Laage2021}, we start by describing their PDE formulation highlighting the gap we address. Our proposed two-step methodology combines a data-driven first step which can be seen as clustering with a heuristic in the second step. The novelty resides in this combination that is made possible thanks to a extension of the formulation in \cite{Laage2021}. We do not claim any contribution to the literature on clustering or on heuristics. Therefore, in Section~\ref{sec:litSolutionApp}, we only provide a brief background on these topics with pointers to work we use in this paper.

\subsection{The Periodic Demand Estimation Problem}

\cite{Laage2021} introduce the PDE formulation which translates the tactical planning process of a freight carrier into a multilevel optimization problem. The objective is to minimize the tactical costs over the tactical planning horizon by finding a good estimate of the periodic demand.

We consider a tactical planning horizon that is decomposed into $T$ periods indexed by $t=1,\ldots,T$. 
The transport network carries a set of commodities $\mathcal{K}$, and each commodity $k$ is described by its origin $o_k$, destination $d_k$ and type $\gamma_k$. We denote $K$ the cardinality of $\mathcal{K}$. 
The demand vector of period $t$ is $\mathbf{y}_t = ( y_{t1}, \dots, y_{tK})^{\top}$, where $y_{tk}$ is the quantity of commodity $k$ to be transported during period $t$. 
Let $\mathbf{Y} \in \mathbb{R}_+^{T \times K}$ be the demand matrix, with $[\mathbf{Y}]_{tk} = y_{tk}$.
Similarly, the periodic demand vector is $\mathbf{y}^{\text{p}} = ( y_1^{\text{p}}, \dots, y_K^{\text{p}})^{\top}$, where $y_k^{\text{p}}$ is the periodic demand for commodity $k$. Let us note $\mathcal{Y}$ the set of feasible values for $\mathbf{y}^{\text{p}}$. In practice, the demand is not known in advance and the demand values are forecasts. 

There are different SND formulations proposed in the literature.
\cite{Laage2021} use the Multicommodity Capacitated Fixed-charge Network Design (MCND) problem \citep{magnanti1984network} for illustration purposes, and we do the same in this paper. 
The generic formulation of a path-based MCND relies on a space-time graph $\mathcal{G} = (\mathcal{N}, \mathcal{A})$ where $\mathcal{N}$ is the set of nodes and $\mathcal{A}$ is the set of arcs. 
A path $p$ is a sequence of arcs in $\mathcal{G}$, and $\mathcal{P}$ denotes the set of paths. Let $\mathcal{P}_k$ denote the set of paths for commodity $k$ such that the source node of the first arc of $p \in \mathcal{P}_k$ is $o_k$ and the sink node of the last arc is $d_k$. We denote $\mathcal{K}_p$ the set of commodities that can use path $p \in \mathcal{P}$, and  $\mathcal{P}^{\text{out}} \subset \mathcal{P}$ the set of paths that correspond to outsourcing options transporting demand in case of insufficient capacity.
Similarly, the set $\mathcal{P}_k^{\text{out}}\subset \mathcal{P}_k$ denotes the outsourcing paths for commodity $k$.

We present below the formulation \textbf{PDE}, where \textbf{MCND} and \textbf{wMCND} constitute the lower levels. The upper level aims at minimizing the total fixed and variable costs over the tactical planning horizon, and $\mathbf{y}^{\text{p}}$ is the decision variable. 
The objective of \textbf{MCND} is to satisfy the periodic demand at minimum cost. There are two types of decision variables: the design variables and the flow variables. The binary design variables $z_p, ~\forall p\in \mathcal{P},$ are equal to 1 if path $p$ is used and 0 otherwise, and the flow variables $x_{pk} \geq 0,  p\in \mathcal{P}_k, ~\forall k\in \mathcal{K}$ are either continuous or integers depending on the type of freight. The third level \textbf{wMCND} aims at satisfying demand at each period of the horizon at minimum cost for a fixed design solution $z$ given by \textbf{MCND}. The variables $x_{tpk}$ equal the flow for commodity $k$ on path $p$ in period $t$.

\begin{footnotesize}
\begin{flalign}
    \textbf{PDE} \quad \min_{\mathbf{y}^{\text{p}}} & \quad C^{\text{PDE}} = \sum_{t=1}^T \left[ \sum_{p \in \mathcal{P}} C_p^{\text{design}} z_p  +  \sum_{k \in \mathcal{K}} \sum_{p \in \mathcal{P}_k \setminus \mathcal{P}_k^{\text{out}}} C_p^{\text{flow}} x_{tpk} +  \sum_{k \in \mathcal{K}} \sum_{p \in \mathcal{P}_k^{\text{out}}}  C_p^{\text{out}} x_{tpk} \right] \label{eq:genericP_obj} \\
    %
    \text{s.t.} & \hspace{0.3cm} \mathbf{y}^{\text{p}} \in \mathcal{Y}, \label{eq:yp_generic} \\
    \textbf{MCND} & \hspace{0.3cm} \min_{z, x} \sum_{p \in \mathcal{P}} C_p^{\text{design}} z_p   +  \sum_{k \in \mathcal{K}} \sum_{p \in \mathcal{P}_k \setminus \mathcal{P}_k^{\text{out}}} C_p^{\text{flow}} x_{pk} +  \sum_{k \in \mathcal{K}} \sum_{p \in \mathcal{P}_k^{\text{out}}}  C_p^{\text{out}} x_{pk} \label{eq:mcnd_objective} \\
    & \hspace{0.3cm} \text{ s.t. } \sum_{p \in \mathcal{P}_k} x_{pk} = y_k^{\text{p}}, \hspace{2.3cm} k \in \mathcal{K}, \label{eq:mcnd_demand}\\ 
    & \hspace{0.8cm} \sum_{k \in \mathcal{K}_p} x_{pk} \leq u_p z_p, \hspace{2cm} p \in \mathcal{P}, \label{eq:mcnd_cap}\\ 
    & \hspace{0.8cm} x_{pk} \geq 0 , \hspace{3.2cm} k \in \mathcal{K}, p \in \mathcal{P}_k,  \label{eq:mcnd_flowPos}\\ 
    & \hspace{0.8cm} z_p \in \{0,1\}, \hspace{2.8cm}  p \in \mathcal{P},  \label{eq:mcnd_design}\\
    \textbf{wMCND} & \hspace{0.8cm} \min_{x_1,\ldots,x_T} \sum_{t=1}^T \left[ \sum_{k \in \mathcal{K}} \sum_{p \in \mathcal{P}_k \setminus \mathcal{P}_k^{\text{out}}} C_p^{\text{flow}} x_{tpk} +  \sum_{k \in \mathcal{K}} \sum_{p \in \mathcal{P}_k^{\text{out}}}  C_p^{\text{out}} x_{tpk} \right] \label{eq:genericP_l2_obj} \\
    & \hspace{0.8cm} \text{ s.t. } \sum_{p \in \mathcal{P}_k} x_{tpk} = y_{tk}, \hspace{1.6cm} t=1,\ldots,T, k \in \mathcal{K}, \label{eq:wmcnd_demand}\\ 
    & \hspace{1.4cm} \sum_{k \in \mathcal{K}_p} x_{tpk} \leq u_p z_p, \hspace{1.3cm} t=1,\ldots,T, p \in \mathcal{P}, \label{eq:wmcnd_cap}\\ 
    & \hspace{1.4cm} x_{tpk} \geq 0 , \hspace{2.6cm} t=1,\ldots,T, k \in \mathcal{K}, p \in \mathcal{P}_k. \label{eq:wmcnd_flowPos}
\end{flalign}
\end{footnotesize}

The objective function~\eqref{eq:genericP_obj} of \textbf{PDE} is the sum of the costs over the tactical planning horizon and depend on the design and flow variables from the lower levels \textbf{MCND} and \textbf{wMCND}. The solution is constrained to $\mathcal{Y}$ which we discuss in more detail below.

The objective function~\eqref{eq:mcnd_objective} of \textbf{MCND} includes three terms. The first term corresponds to the fixed design cost $C_p^{\text{design}}\geq0$ and accounts for the paths built to transport demand. The second term is the variable flow cost $C_p^{\text{flow}}\geq0$ for satisfied demand, and the third term with $C_p^{\text{out}}\geq0$ corresponds to the flow cost of outsourced demand. The third-level objective function~\eqref{eq:genericP_l2_obj} contains the last two terms, that is, the variable flow cost and the outsourcing cost at each period. 
Furthermore, the two lower levels contain two sets of the same type of constraints. Namely, constraints~\eqref{eq:mcnd_demand} and \eqref{eq:wmcnd_demand} ensure that, respectively, the periodic demand and the demand at each period are satisfied for each commodity. Constraints~\eqref{eq:mcnd_cap} and~\eqref{eq:wmcnd_cap} enforce flows on selected paths only, and that, respectively, the flow from the periodic demand and the flow at each period do not exceed the path capacity $u_p$. 

In order to simplify the description of our methodology, we introduce a compact form of the \textbf{PDE} formulation:
\begin{small}
\begin{align}
    \textbf{PDE} \quad \min_{\mathbf{y}^{\text{p}}} & \quad C^{\text{PDE}}(\mathbf{y}^{\text{p}}, \mathbf{z}, \mathbf{x}, \mathbf{x_1},\ldots,\mathbf{x_T}) \label{eq:yp_genericObj}\\
    \text{s.t.} & \hspace{0.3cm} \mathbf{y}^{\text{p}} \in \mathcal{Y}, \label{eq:yp_generic2}\\
    & \hspace{0.3cm} (\mathbf{z}, \mathbf{x}, \mathbf{x_1},\ldots,\mathbf{x_T}) \in \argmin_{\mathbf{z}', \mathbf{x}', \mathbf{x_1}',\ldots,\mathbf{x_T}'} \textbf{MCND-wMCND}(\mathbf{y}^{\text{p}}, \mathbf{z}', \mathbf{x}', \mathbf{x_1}',\ldots,\mathbf{x_T}') \label{eq:yp_generic3}
\end{align}
\end{small}

Under a weak assumption, the formulation \textbf{PDE} does not belong to the class of multilevel problems, despite its appearance. For $(z^*, x^*)$ an optimal solution of \textbf{MCND}, if $z^*$ is feasible for \textbf{wMCND}, then \textbf{MCND}-\textbf{wMCND} can be solved sequentially for a fixed $\mathbf{y}^{\text{p}}$. We call this the \textit{sequential property} \citep[Claim~1 in][]{Laage2021}. 
By always allowing outsourcing, i.e., either fixing $z_p = 1$ or $C_p^{\text{design}} = 0$ for $p \in \mathcal{P}^{\text{out}}$, the property holds. 

The periodic demand is defined as a mapping of the demand values at each week of the tactical horizon. The set of feasible periodic demand $\mathcal{Y}$ is then a set of feasible mappings $h$, where $h$ is defined as 
\begin{equation}
    \begin{split}
    h \colon \quad \mathbb{R}_+^{T \times K} & \to \mathbb{R}_+^{K}  \\
    \mathbf{Y} &\mapsto \mathbf{y}^{\text{p}} = h(\mathbf{Y}).
    \end{split}
    \label{eq:h_mapping}
\end{equation}
Even if the sequential property is satisfied, the problem is challenging to solve because lower levels are large-scale, non-convex, non-differentiable and combinatorial problems. \cite{Laage2021} restrict $\mathcal{Y}$ to the following four mappings,
\begin{equation}
    \mathbf{y}_{\text{max}}^{\text{p}} = h_1(\mathbf{y}_1, \ldots, \mathbf{y}_T)= \max_{t=1,\ldots,T} \mathbf{y}_t,
    \label{eq:max_periodic}
\end{equation}
\begin{equation}
    \mathbf{y}_{\text{mean}}^{\text{p}} =  h_2(\mathbf{y}_1, \ldots, \mathbf{y}_T) =  \frac{1}{T} \sum_{t=1}^T \mathbf{y}_t,
    \label{eq:mean_periodic}
\end{equation}
\begin{equation}
    \mathbf{y}_{\text{q2}}^{\text{p}} =  h_3(\mathbf{y}_1, \ldots, \mathbf{y}_T) =  Q_2(\mathbf{y}_{t}, t = 1 \ldots, T),
    \label{eq:q2_periodic}
\end{equation}
\begin{equation}
    \mathbf{y}_{\text{q3}}^{\text{p}} =  h_4(\mathbf{y}_1, \ldots, \mathbf{y}_T) =  Q_3(\mathbf{y}_{t}, t = 1 \ldots, T),
    \label{eq:q3_periodic}
\end{equation}
representing the maximum, mean, second quartile $Q_2$ and third quartile $Q_3$, respectively. This restriction of the solution space makes it possible to solve \textbf{PDE} by enumeration and \cite{Laage2021} report significant cost reductions compared to standard practice (i.e., using $\mathbf{y}_{\text{mean}}^{\text{p}}$). However, the actual solution space is not restricted to only a few mappings and hence even higher quality solutions can exist. In this work we propose a new formulation based on a continuous extension of the solution space. In turn, this leads to a problem that cannot be solved by enumeration and we devise a two-step heuristic combining a data-driven clustering heuristic with either black-box optimization algorithms or metaheuristics.

\subsection{Work Related to the Solution Approach} \label{sec:litSolutionApp}

In this section we provide a brief and high-level background on three topics related to our work, namely, clustering, metaheuristics and black-box optimization.

Clustering is an unsupervised learning method consisting in grouping data points \citep[see, e.g.,][for a comprensive survey]{xu2015comprehensive}. More precisely, given a set of multi-dimensional data points, clustering algorithms assign a specific group to each data point with the objective to group points with similar properties. Popular algorithms include $k$-means \citep{macqueen1967some} and Expectation-Maximization clustering using Gaussian Mixture Models \citep{rasmussen1999infinite}. 
To preserve the interpretability of solutions in this work, we propose tailored clustering heuristics based on the data analysis of the problem structure, rather than using unsupervised learning algorithms.

In the second step of our proposed two-step heuristic, we can use a metaheuristic or a black-box optimization algorithm to solve \textbf{PDE}.
Both classes of methods search for high-quality solutions by iterating over a sequence of candidate solutions. In the context of our work, it is important to consider the large size, non-convexity and non-differentiability of \textbf{PDE}'s lower levels. Indeed, for each candidate solution, we need to solve the lower levels \textbf{MCND}-\textbf{wMCND} to evaluate $C^{\text{PDE}}$~\eqref{eq:genericP_obj} which, in turn, directs the search. 

Metaheuristics can coarsely be divided in two types \citep[e.g.,][]{talbi2009metaheuristics}: local search methods iterating over one solution and population methods working with a population of solutions. Local search methods are better suited for our problem as each solution can be expensive to evaluate. There is an extensive literature on this topic with many types of metaheuristics: Tabu Search \citep{glover1986}, Simulated Annealing \citep{kirkpatrick1983optimization, cerny1985}, Variable Neighborhood Search \citep{mladenovic1997}, Greedy Randomized Adaptive Search Procedure \citep{feo1989, feo1995}, to name a few. They differ by the set of solutions visited at each iteration, the use of memory during the search and the diversification and intensification procedures to explore the feasible space. 

Black-box optimization methods are also extensively studied in the literature \citep[e.g.,][]{audet2017derivative}. 
We take a particular interest in NOMAD \citep{le2011algorithm}, a state-of-the-art black-box optimization solver  implementing the Nonsmooth Optimization by Mesh Adaptive Direct Search \citep{abramson2004filter,abramson2009mesh}. It is openly accessible and designed to be stable and to address real problems, therefore appropriate for our framework. As most black-box optimization algorithms, NOMAD performs best when there is a relatively low number of variables, less than 50 according to \cite{le2011algorithm}. 

In terms of methodological contribution, the novelty in our approach consists in the formulation and the first, data-driven step of the two-step heuristic. In the second step, we use existing results. More precisely, we propose metaheuristics as baselines by using well-known ideas from that literature, and we compare the performance with NOMAD.

\section{Methodology}\label{section:methodology}

In this section we start by introducing a formulation of \textbf{PDE} where the periodic demand is defined as a deviation from the average of the demand forecasts. We then illustrate the importance of the problem and highlight differences with respect to \cite{Laage2021} by means of a small example (Section~\ref{sec:IllEx}). In Section~\ref{section:clustering}, we describe the first step of our methodology: clustering heuristics that restrict the solution space $\mathcal{Y}$ in a data-driven way. Finally, in Section~\ref{section:metaheuristics}, we present the local search metaheuristics used as baselines in the second step to solve the PDE problem. 

\subsection{A Formulation Based on Deviation Coefficients}
\label{section:model}

The key idea behind the formulation we introduce in this paper is the extension of the feasible space, compared to the discrete set of solutions proposed in \cite{Laage2021}. To ensure computational tractability, we then restrict the space in a data-driven way with the aim to be less restrictive than \cite{Laage2021}. For this purpose, we
define the periodic demand as a deviation from $\mathbf{y}_{\text{mean}}^{\text{p}}$,
\begin{equation}
    \mathbf{y}^{\text{p}} =  \boldsymbol{\alpha} \odot \mathbf{y}_{\text{mean}}^{\text{p}}, 
\end{equation}
where $\odot$ designates the element-wise multiplication. The vector $\boldsymbol{\alpha} \in \mathbb{R}_+^K$ of what we refer to as \emph{deviation coefficients} is defined such that, for each commodity $k$, 
$y_k^{\text{p}} = [\boldsymbol{\alpha}]_k y_{\text{mean},k}^{\text{p}}$. We assume that the periodic demand values are bounded by the minimum and the maximum demand values. Hence, the solution space for the periodic demand, expressed over $\boldsymbol{\alpha}$ is
\begin{equation}
    \mathcal{Y} = \{\boldsymbol{\alpha}, \boldsymbol{\alpha}_{\text{min}} \leq \boldsymbol{\alpha} \leq \boldsymbol{\alpha}_{\text{max}}\}, 
\end{equation}
where $\boldsymbol{\alpha}_{\text{max}}$ is defined such that $\mathbf{y}_{\text{max}}^{\text{p}} =  \boldsymbol{\alpha}_{\text{max}} \odot \mathbf{y}_{\text{mean}}^{\text{p}}$, with $\mathbf{y}_{\text{max}}^{\text{p}} = \max_{t=1,\ldots,T} \{\mathbf{y}_t\}$ and $\boldsymbol{\alpha}_{\text{min}}$ is defined such that $\mathbf{y}_{\text{min}}^{\text{p}} =  \boldsymbol{\alpha}_{\text{min}} \odot \mathbf{y}_{\text{mean}}^{\text{p}}$, with $\mathbf{y}_{\text{min}}^{\text{p}} = \min_{t=1,\ldots,T} \{\mathbf{y}_t\}$. 

We can now introduce the relaxed formulation \textbf{ePDE} where the upper level variables are the deviation coefficient vector $\boldsymbol{\alpha}$
\begin{small}
\begin{align}
    \textbf{ePDE} \quad \min_{\boldsymbol{\alpha}} & \quad C^{\text{PDE}}(\mathbf{y}^{\text{p}}, \mathbf{z}, \mathbf{x}, \mathbf{x_1},\ldots,\mathbf{x_T}) \label{eq:new_formul1a}\\
    \text{s.t.} & \hspace{0.3cm} \mathbf{y}^{\text{p}} =  \boldsymbol{\alpha} \odot \mathbf{y}_{\text{mean}}^{\text{p}},\label{eq:new_formul2a} \\
    & \hspace{0.3cm} \boldsymbol{\alpha} \leq \boldsymbol{\alpha}_{\text{max}} \label{eq:new_formul3a}\\
    & \hspace{0.3cm} \boldsymbol{\alpha} \geq \boldsymbol{\alpha}_{\text{min}} \label{eq:new_formul4a}\\
    & \hspace{0.3cm} (\mathbf{z}, \mathbf{x}, \mathbf{x_1},\ldots,\mathbf{x_T}) \in \argmin_{\mathbf{z}', \mathbf{x}', \mathbf{x_1}',\ldots,\mathbf{x_T}'} \textbf{MCND-wMCND}(\mathbf{y}^{\text{p}}, \mathbf{z}', \mathbf{x}', \mathbf{x_1}',\ldots,\mathbf{x_T}') \label{eq:new_formul5a}
\end{align}
\end{small}
where Constraints~\eqref{eq:new_formul3a} and ~\eqref{eq:new_formul4a} respectively impose upper and lower bounds on $\boldsymbol{\alpha}$. 

A few observations are in order. First, \eqref{eq:genericP_obj}-\eqref{eq:wmcnd_flowPos} with \eqref{eq:max_periodic}-\eqref{eq:q3_periodic} is a special case of \eqref{eq:new_formul1a}-\eqref{eq:new_formul5a}. Second, the dimension of $\boldsymbol{\alpha}$ is the number of commodities $K$ and hence can be large, in our application $K=170$. Given the difficulty of problems \eqref{eq:new_formul5a} this implies that it may not be possible to solve \textbf{ePDE} for the full dimension of $\boldsymbol{\alpha}$. This is the motivation for our clustering heuristic discussed in Section~\ref{section:clustering}. Third, the solution is easy to interpret. As the name entails, for each commodity, the deviation coefficient represents a gap between the periodic and average demand.

\subsection{Illustrative Example} \label{sec:IllEx}
\begin{figure}[!htbp]  
\centering
\begin{tikzpicture}
    \node[shape=circle,draw=black, fill = gray!50, scale = 0.8] (O1) at (1,2) {$O_1$};
    \node[shape=circle,draw=black, fill = gray!20, scale = 0.8] (T1) at (10, 2) {$D_1$};
    
   \path [bend right = 40, ->] (O1) edge node[below] {$u_1 = 2$, $C^{\text{design}}_1 = 10$, $C^{\text{flow}}_1 = 5$} (T1);
   \path [bend right = 10, dashdotted, ->, line width=0.4mm] (O1) edge node[below] {$u_2 = 1$, $C^{\text{design}}_2 = 10$, $C^{\text{flow}}_2 = 10$} (T1);
   \path [bend left = 10, dotted, ->, line width=0.4mm] (O1) edge node[above] {$u_3 = 2$, $C^{\text{design}}_3 = 20$, $C^{\text{flow}}_3 = 10$} (T1);
   \path [bend left = 40, densely dashed, ->, line width=0.4mm] (O1) edge node[above] {$u_4 = \infty$, $C^{\text{design}}_4 = 0$, $C^{\text{out}}_4 = 50$ } (T1);

\end{tikzpicture}
\caption{Illustration of the PDE problem on a small network}
\label{fig:network_illustration}
\end{figure}
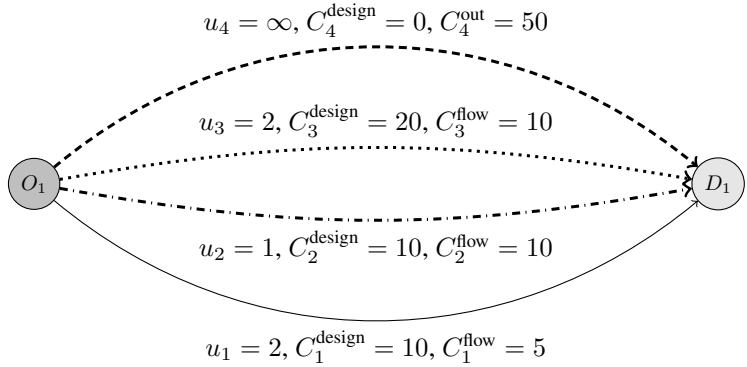

In this section we introduce a small example to illustrate the importance of the PDE problem as well as the potential cost reduction by solving \textbf{ePDE}~\eqref{eq:new_formul1a}-\eqref{eq:new_formul5a} compared to the formulation in \cite{Laage2021}.
We consider a network transporting one commodity, from its origin $O_1$ to its destination $D_1$ over a tactical horizon of six periods. Figure~\ref{fig:network_illustration} presents the four different possible paths $p$ for the commodity, along with their capacity $u_p$, design cost $C_p^{\text{design}}$ and flow cost $C_p^{\text{flow}}$. Path 4 corresponds to the outsourcing path and we denote $C_p^{\text{out}}$ its flow cost. It is always built in the design (the design cost is zero), it is expensive to use but does not have a capacity constraint.

We assume that the time-series demand forecasts $\hat{\mathbf{Y}}$ and the observed demand $\mathbf{Y}$ of the commodity at the six periods are
\begin{equation}
   \hat{\mathbf{Y}} = \begin{pmatrix}  4 \\ 2 \\ 1 \\ 0 \\ 1 \\ 4 \end{pmatrix} \quad \text{and} \quad \mathbf{Y} = \begin{pmatrix}  1 \\  2\\ 1\\ 1 \\ 3\\ 3\\ \end{pmatrix}.
\end{equation}

We report
the results in Table~\ref{table:illustration} where $C^{\text{design}}$ is the design costs of \textbf{MCND}, i.e., the first term of Equation~\eqref{eq:mcnd_objective}. Moreover, $C^{\text{flow}}_{\textbf{MCND}}$ is the sum of the flow and outsourcing costs of \textbf{MCND}, i.e., the second and third term of Equation~\eqref{eq:mcnd_objective}, and $C_{\textbf{wMCND}}$ is the value of the objective function of \textbf{wMCND}~\eqref{eq:genericP_l2_obj}. The second last column reports the tactical costs~\eqref{eq:genericP_obj}, for a horizon of six periods, it is equal to $6 C^{\text{design}} + C_{\textbf{wMCND}}$.  In the last column, we report the actual tactical costs $C_{\text{act}}^{\text{PDE}}$ computed after the demand realizations. 

\begin{table}[!htbp]
\centering
\begin{tabular}{l|llllll}
Periodic Demand & Paths built & $C^{\text{design}}$ & $C^{\text{flow}}_{\textbf{MCND}}$ & $C_{\textbf{wMCND}}$ & $C^{\text{PDE}}$ & $C_{\text{act}}^{\text{PDE}}$\\ \hline
$\hat{y}_{\text{mean}}^{\text{p}} = 2$ & Path 1        & 10 & 10 & 240 & 300 & 205 \\
$\hat{y}_{\text{q2}}^{\text{p}} = 2$   & Path 1        & 10 & 10 & 240 & 300 & 205 \\
$\hat{y}_{\text{q3}}^{\text{p}} = 4$   & Paths 1,2,3   & 40 & 30 & 80 & 320 & 305 \\
$\hat{y}_{\text{max}}^{\text{p}} = 4$  & Paths 1,2,3   & 40 & 30 & 80 & 320 & 305 \\ \hdashline
$\alpha = 1.5$, $\hat{y}^{\text{p}} = 3$ & Paths 1, 2 & 20 & 20 & 160 & 280 & 185
\end{tabular}
\caption{Results for the illustrative example (the second and third quartiles are rounded up to the nearest integer)}
\label{table:illustration}
\end{table}

The first four rows of the table reports the demand for each of the mappings considered in \cite{Laage2021}. The last row corresponds to a deviation coefficient $\alpha = 1.5$, i.e., $\hat{y}^{\text{p}} =  \alpha \hat{y}_{\text{mean}}^{\text{p}} = 3$.

The results show that the periodic demand using $\alpha=1.5$

leads to the lowest cost. As expected, relaxing the solution space can be valuable. 
Moreover, the value of the deviation coefficient $\alpha$ allows to easily interpret the solution: to minimize the costs, the periodic demand should be 50\% higher than the average forecast.

\subsection{Clustering to Restrict the Feasible Space}
\label{section:clustering}

In this section we propose a data-driven clustering heuristic, which consists in creating clusters of commodities to reduce the dimension of the feasible space. The objective of clustering is to reduce the number of variables while keeping high-quality solutions. 
More precisely, the deviation coefficients of commodities grouped into a same cluster are assumed to be equal.
It corresponds to adding equality constraints to \eqref{eq:new_formul1a}-\eqref{eq:new_formul5a} (heuristic variable fixing) and hence restricting the solution space.

Intuitively, the set of clusters should depend on the application under consideration. For instance, increased tactical costs can be due to the outsourcing of commodities. This happens when the periodic demand used to define the plan is low and the network does not have sufficient capacity. A cluster could gather the commodities with high risks of being outsourced, so that \textbf{ePDE} assigns them a large deviation coefficient. Before describing two clustering heuristics that use data on the network and the demand (Subsections~\ref{section:variance clusters} and~\ref{subsection:resource clustering}), we introduce some basic notation.

We denote $\mathcal{C} = \{ C_1, \ldots, C_{n_C} \}$ the set of $n_C$ clusters. It is a partition of $\mathcal{K}$ such that for $ C_i \in \mathcal{C}$ with $|C_i| = c$ and $C_i = \{ k_1, .., k_c \}$, the coefficients $\alpha_{k'}$ for the commodities $k' \in C_i$ are equal. 
This leads to the formulation \textbf{cPDE} with clusters:
\begin{small}
\begin{align}
    \textbf{cPDE} \quad \min_{\boldsymbol{\alpha}} & \quad C^{\text{PDE}}(\mathbf{y}^{\text{p}}, \mathbf{z}, \mathbf{x}, \mathbf{x_1},\ldots,\mathbf{x_T}) \label{eq:new3_formul1}\\
    \text{s.t.} & \hspace{0.3cm} \mathbf{y}^{\text{p}} =  \boldsymbol{\alpha} \odot \mathbf{y}_{\text{mean}}^{\text{p}},\label{eq:new3_formul2} \\
    & \hspace{0.3cm} \boldsymbol{\alpha} \leq \boldsymbol{\alpha}_{\text{max}} \label{eq:new3_formul3a}\\
    & \hspace{0.3cm} \boldsymbol{\alpha} \geq \boldsymbol{\alpha}_{\text{min}} \label{eq:new3_formul4a}\\
    & \hspace{0.3cm} [\boldsymbol{\alpha}]_{k_i} = [\boldsymbol{\alpha}]_{k_j}, \hspace{1.3cm}   k_i, k_j \in C, C \in \mathcal{C}, \label{eq:new3_formul4b} \\
    & \hspace{0.3cm} (\mathbf{z}, \mathbf{x}, \mathbf{x_1},\ldots,\mathbf{x_T}) \in \argmin_{\mathbf{z}', \mathbf{x}', \mathbf{x_1}',\ldots,\mathbf{x_T}'} \textbf{MCND-wMCND}(\mathbf{y}^{\text{p}}, \mathbf{z}', \mathbf{x}', \mathbf{x_1}',\ldots,\mathbf{x_T}'). \label{eq:new3_formul5}
\end{align}
\end{small}

One simple example is to constrain all variables to be equal, i.e., forming only one cluster. This is equivalent to having a single decision variable $\alpha$ and the formulation reduces to
\begin{small}
\begin{align}
    \quad \min_{\alpha} & \quad C^{\text{PDE}}(\mathbf{y}^{\text{p}}, \mathbf{z}, \mathbf{x}, \mathbf{x_1},\ldots,\mathbf{x_T}) \label{eq:new2_formul1}\\
    \text{s.t.} & \hspace{0.3cm} \mathbf{y}^{\text{p}} =  \alpha \mathbf{y}_{\text{mean}}^{\text{p}},\label{eq:new2_formul2} \\
    & \hspace{0.3cm} \alpha \leq \max_k \boldsymbol{\alpha}_{\text{max},k} \label{eq:new2_formul3}\\
    & \hspace{0.3cm} \alpha \geq \min_k \boldsymbol{\alpha}_{\text{min},k} \label{eq:new2_formul4}\\
    & \hspace{0.3cm} (\mathbf{z}, \mathbf{x}, \mathbf{x_1},\ldots,\mathbf{x_T}) \in \argmin_{\mathbf{z}', \mathbf{x}', \mathbf{x_1}',\ldots,\mathbf{x_T}'} \textbf{MCND-wMCND}(\mathbf{y}^{\text{p}}, \mathbf{z}', \mathbf{x}', \mathbf{x_1}',\ldots,\mathbf{x_T}'). \label{eq:new2_formul5}
\end{align}
\end{small}

\subsubsection{Variance-based Clustering}
\label{section:variance clusters}

We propose a first clustering based on the coefficient of variation, i.e., the standard deviation scaled to the average demand to prevent the increase of costs due to outsourcing. Indeed, let us assume for instance that we take $\boldsymbol{\alpha} = \mathbf{1}_K$ the $K$-vector of ones, i.e., $\mathbf{y}^{\text{p}} = \mathbf{y}_{\text{mean}}^{\text{p}}$ to build the service network.
If the demand forecasts for a commodity have a large variance over the planning horizon, then at certain periods the demand realizations are large compared to the periodic demand. This results in outsourcing of the commodity, if the network does not have sufficient capacity. 

We denote $\sigma_k^{\text{s}}$ the coefficient of variation for a commodity $k$, such that
\begin{equation}
    \sigma_k^{\text{s}} = \frac{\sqrt{\frac{1}{T} \sum_{t=1}^T (y_{tk} - y_{\text{mean, k}}^{\text{p}})^2 }}{y_{\text{mean, k}}^{\text{p}}}, \quad k = 1, \ldots, K.
\end{equation}

When $K$ is small, we can create the clusters from analyzing the distribution of $\sigma^{\text{s}}$ over the commodities. More generally, the $n_C$ clusters are created by decomposing the set $\{ \sigma_k^{\text{s}} , k=1,\ldots,K\}$ into $n_C$ intervals. Let $Q_e$ denote the e-th percentile of $\{ \sigma_k^{\text{s}} , k=1,\ldots,K\}$. For instance we can create $n_C = 2$ clusters, where the first cluster contains all $\sigma_k^{\text{s}} \leq Q_{0.25}$ and the second cluster contains all $\sigma_k^{\text{s}} > Q_{0.25}$. In our application, we consider $n_C = 5$ clusters and describe them in Section~\ref{section:results}.

\subsubsection{Service-based Clustering}
\label{subsection:resource clustering}

We propose a second approach built from the analysis of the shared services or resources in the network. Increasing the periodic demand for commodities that share resources can lead to an increased dedicated capacity in the design. This, in turn, might help to handle differences between forecasts and observed demand.
We propose the following steps to identify the set of clusters.

\begin{enumerate}
    \item We start from the common practice where the periodic demand is the average of the demand forecasts, i.e., we compute $C^{\text{PDE}} (\boldsymbol{\alpha}_{\text{mean}})$ by solving $\textbf{ePDE}$ where $\boldsymbol{\alpha}_{\text{mean}} = \mathbf{1}_K$.
    \item For each commodity $k$, we identify the commodities that share at least one service (a train for instance) on their route with $k$. We call a \textit{group} the set formed by $k$ and the related commodities. 
    \item The first cluster is the largest group over all commodities. 
    \item The second cluster is the second largest group of cardinality higher than 1 that has no commodities in common with the first cluster. We iterate until there is no more group satisfying the criteria. 
    \item The remaining commodities not yet assigned to a cluster are gathered to form the last cluster. 
\end{enumerate}

\paragraph{Relaxed Service-based Clustering}
A similar alternative consists in modifying the first step.
Instead of computing $C^{\text{PDE}} (\boldsymbol{\alpha}_{\text{mean}})$ from the formulation \textbf{ePDE}, we relax the capacity constraints in the lower levels and solve this new formulation with $\boldsymbol{\alpha}_{\text{mean}}$. This aims to identify the best paths to transport each commodity, and analyze the bottlenecks in the network. We then follow the steps 2 to 5 to obtain the clusters.

\subsection{Metaheuristics}
\label{section:metaheuristics}

In this section we focus on local search methods to solve \textbf{ePDE}:
the first-level constraints help to define the set of solutions to visit and we use the sequential property of \textbf{ePDE} at each iteration, when a new vector $\boldsymbol{\alpha}$ is selected. We fix the upper-level variable to $\boldsymbol{\alpha}$, then solve \textbf{MCND-wMCND} sequentially and finally use the solution to compute $C^{\text{PDE}}$, whose value indicates where to continue the search. We present next two tailored local search metaheuristics designed to solve \textbf{ePDE}. They serve as baselines to a state-of-the-art black-box optimization solver. 
Both require a feasible space of vectors that can be visited, a neighborhood $N$ indicating the movements allowed in the feasible space and a stopping criterion. 
We present below the pseudo-codes and related parameters. 

\subsubsection{Neighborhood Search}

We name Neighborhood Search (NS) the first metaheuristic, which is a simple local search method. 
At each iteration, it explores the neighborhood around a given solution and stops when the value of $C^{\text{PDE}}$ is not improved. We present in Algorithm~\ref{algo:descent} the pseudo-code for NS and in Algorithm~\ref{algo:neighbors} the pseudo-code for the neighborhood definition.

\begin{algorithm}[H]
\caption{Neighborhood Search (NS)}
\label{algo:descent}
\begin{algorithmic}[1]
\State \textbf{Input}: Initial solution $\boldsymbol{\alpha}$
\State Stop = False
\State \textbf{while} not Stop \textbf{do}
\State \hspace{1cm} Define $N(\boldsymbol{\alpha})$ the neighborhood around the solution $\boldsymbol{\alpha}$ with Algorithm~\ref{algo:neighbors}
\State \hspace{1cm} Find the solution $\boldsymbol{\alpha'}$ minimizing $C^{\text{PDE}}$ in $N(\boldsymbol{\alpha})$
\State \hspace{1cm} \textbf{If} $C^{\text{PDE}}(\boldsymbol{\alpha'}) < C^{\text{PDE}}(\boldsymbol{\alpha})$
\State \hspace{2cm} \textbf{Update} $\boldsymbol{\alpha} := \boldsymbol{\alpha'}$
\State \hspace{1cm} \textbf{Else} Stop := True
\end{algorithmic}
\end{algorithm}

The elements of the neighborhood $N(\boldsymbol{\alpha})$ for a solution $\boldsymbol{\alpha}$ are randomly generated following a Normal distribution. We denote $V$ the number of neighbors, i.e., $V = | N(\boldsymbol{\alpha})|$, $\beta$ a distance parameter and $\mathbf{I}_{K}$ the identity matrix of size $K$. Step 6 in Algorithm~\ref{algo:neighbors} ensures that the potential solution vectors satisfy Constraints~\eqref{eq:new_formul3a} and \eqref{eq:new_formul4a}. 

\begin{algorithm}[H]
\caption{Neighborhood Definition}
\label{algo:neighbors}
\begin{algorithmic}[1]
\State \textbf{Input}: $\boldsymbol{\alpha}$
\State \textbf{Parameters}:  $\beta$, $V$
\State \textbf{Initialization}  $N(\boldsymbol{\alpha}) = \emptyset$
\State \textbf{while} $|N(\boldsymbol{\alpha})| < V$ \textbf{do}
\State \hspace{1cm}  Generate random $\boldsymbol{\alpha'}$ from the normal distribution $\mathcal{N}(\boldsymbol{\alpha}, \beta \mathbf{I}_{K})$
\State \hspace{1cm}  $\boldsymbol{\alpha'} := \max(\boldsymbol{\alpha'}, \boldsymbol{\alpha}_{\text{min}}$), $\boldsymbol{\alpha'} := \min(\boldsymbol{\alpha'}, \boldsymbol{\alpha}_{\text{max}})$
\State \hspace{1cm} \textbf{Update} $N(\boldsymbol{\alpha}) := N(\boldsymbol{\alpha}) \cup \{ \boldsymbol{\alpha'} \}$
\State \textbf{Output} $N(\boldsymbol{\alpha})$
\end{algorithmic}
\end{algorithm}

\subsubsection{Neighborhood Search with Diversification and Intensification}

The well-known limitation of simple algorithms such as NS is their propensity to get stuck at the first local minimum. To address this challenge, we use a metaheuristic allowing at each iteration to move to a neighbor that might not improve the best found objective function. We call it the Neighborhood Search with Diversification and Intensification (NSDI). The neighborhood $N$ is defined as in NS, described in Algorithm~\ref{algo:neighbors}. We add a diversification and intensification procedure, which consists in updating the parameters $\beta$ and $V$ at each iteration whether a better solution was found or not. This metaheuristic, presented next in Algorithm~\ref{algo:improved_greedy}, stops when a maximum number of iterations without improvements $M$ is reached.

\begin{algorithm}[H]
\caption{Neighborhood Search with Diversification and Intensification (NSDI)}
\label{algo:improved_greedy}
\begin{algorithmic}[1]
\State \textbf{Input}: Initial solution $\boldsymbol{\alpha}$ and initial best known solution $\boldsymbol{\alpha}^*$
\State \textbf{Parameters} $\beta$, $V$, $M$, $v^+$, $b^-, b^+$
\State Stop = 0
\State \textbf{while} Stop < $M$ \textbf{do}
\State \hspace{1cm} Define $N(\boldsymbol{\alpha})$ the neighborhood around the solution $\boldsymbol{\alpha}$
\State \hspace{1cm} Find the solution $\boldsymbol{\alpha'}$ which minimizes $C^{\text{PDE}}$ in $N(\boldsymbol{\alpha})$
\State \hspace{1cm} \textbf{If} $C^{\text{PDE}}(\boldsymbol{\alpha'}) < C^{\text{PDE}}(\boldsymbol{\alpha}^*)$
\State \hspace{2cm} \textbf{Update} $\boldsymbol{\alpha}^* := \boldsymbol{\alpha'}$
\State \hspace{2cm} \textbf{Update} Stop:= 0
\State \hspace{2cm} \textbf{Intensification} $\beta := b^- \beta$
\State \hspace{1cm} \textbf{Else} 
\State \hspace{2cm} \textbf{Update} Stop := Stop + 1
\State \hspace{2cm} \textbf{Diversification} $\beta := b^+ \beta$, $V := v^+ V$
\State \hspace{1cm} \textbf{Update} $\boldsymbol{\alpha} := \boldsymbol{\alpha'}$
\end{algorithmic}
\end{algorithm}

In Algorithm~\ref{algo:improved_greedy}, at each iteration, the current solution is updated by the best solution in its neighborhood, even if it does not improve the current value of objective function. Step 10 and 13 perform the diversification and intensification procedure. When the best found solution is improved, we intensify the search in this direction: $b^- < 1$ and generate closer neighbors. On the contrary, when the best found solution is not improved, we look for neighbors that are further away, with $v^+ > 1$ and $b^+ > 1$. 

By definition, NSDI might be stuck between 2 solutions, going back and forth from one to the other. The Tabu Search \citep{glover1986} avoids this problem by defining a tabu list containing solutions that cannot be visited for a given number of iterations. Here, the randomness in the definition of $N$ ensures that this is unlikely to happen and we therefore do not define a tabu list. Indeed, for the same current solution, at different iterations, the set $N$ is likely to be different and the algorithm will move to a different neighbor. The degree of randomness is controlled through the number of neighbors generated and the parameters of the Normal distribution. 

Both metaheuristics NS and NSDI require an initial solution. Since $\boldsymbol{\alpha}$ represents the deviation from the mean, a natural choice for the  initial solution is $\boldsymbol{\alpha} = \mathbf{1}_{K}$, the $K$-vector of ones.

\section{Application}
\label{section:application}

We consider the same application as \cite{Laage2021}. We briefly describe the tactical planning problem here and refer to \cite{Laage2021} for more details on the data.  The application concerns the planning of the intermodal network of the Canadian National Railway Company (CN). It is composed of 24 main terminals and 133 origin-destination pairs. A commodity is defined by an origin, a destination and a type of container. There are two main types of containers, the 40-feet and the 53-feet long, resulting in $K = 170$ commodities. 
The tactical period is a week and a tactical horizon lasts $T = 10$ weeks. 

CN faces a specific MCND problem for its intermodal network, referred to as the \textit{block planning problem}. A block is a consolidation of railcars that move together between a given OD pair, where containers loaded on the railcars share the same OD.
\cite{Morganti2019} introduce a path-based formulation, \textbf{BP}, for this problem where the periodic demand and the schedule of intermodal trains are given as inputs. \textbf{BP} is based on a space-time graph which contains 28,854 arcs and 15,269 nodes. A block is a path in the graph and $\mathcal{B}$ denotes the set of blocks, here $|\mathcal{B}|=2208$. The set $\mathcal{B}$ hence corresponds to the set of paths $\mathcal{P}$ in \textbf{MCND}. It contains a subset of outsourcing paths denoted $\mathcal{B}^{\text{artif}}$, also referred to as \textit{artificial blocks} whose role is to transport demand exceeding capacity.

We now describe \textbf{ePDE} for our application. We refer to \cite{Morganti2019} and \cite{Laage2021} for more details on the lower levels of the formulation, as they correspond to \textbf{BP} and its weekly formulation, \textbf{wBP}. There are three categories of decision variables at the second level. First, the design variables $z_b, b \in \mathcal{B}$ where $z_b = 1$ if block $b$ is built. Second, $x_{bk}$ are integer flow variables corresponding to the number of containers for commodity $k$ transported on block $b$. Third, to write the capacity constraints we need the auxiliary variables for the number of 40-feet $v_b^{40}$ and 53-feet $v_b^{53}$ double-stack platforms on block $b$. At the third level, we introduce flow variables and auxiliary platform variables for each week $t$, $x_{tbk}, v_{tb}^{40}, v_{tb}^{53}, t \in \mathcal{T},  k \in \mathcal{K},  b \in \mathcal{B}$.

\begin{footnotesize}
\begin{flalign}
    \textbf{ePDE} \min_{\boldsymbol{\alpha}} & \hspace{0.02cm} C^{\text{PDE}} = \sum_{t=1}^T \left[ \sum_{b \in \mathcal{B} \setminus \mathcal{B}^{\text{artif}} } C_b^{\text{design}} z_b + \sum_{k \in \mathcal{K}} \sum_{b \in \mathcal{B}_k \setminus \mathcal{B}_k^{\text{artif}} } C_{tbk}^{\text{flow}} x_{tbk} + \sum_{k \in \mathcal{K}} \sum_{b \in \mathcal{B}_k^{\text{artif}} } C_{tbk}^{\text{out}} x_{tbk} \right] \label{eq:obj_pde} \\
    %
    \text{s.t.} & \hspace{0.3cm} \mathbf{y}^{\text{p}} =  \boldsymbol{\alpha} \odot \mathbf{y}_{\text{mean}}^{\text{p}} , \label{eq:pde_bp_2}\\
    %
    & \hspace{0.3cm} \boldsymbol{\alpha} \leq \boldsymbol{\alpha}_{\text{max}} \label{eq:new_formul3}\\
    & \hspace{0.3cm} \boldsymbol{\alpha} \geq \boldsymbol{\alpha}_{\text{min}} \label{eq:new_formul4}\\
    \textbf{BP} & \hspace{0.05cm} \min_{x, z} \sum_{b \in \mathcal{B} \setminus \mathcal{B}^{\text{artif}} } C_b^{\text{design}} z_b + \sum_{k \in \mathcal{K}} \sum_{b \in \mathcal{B}_k \setminus \mathcal{B}_k^{\text{artif}} } C_{bk}^{\text{flow}} x_{bk} + \sum_{k \in \mathcal{K}} \sum_{b \in \mathcal{B}_k^{\text{artif}} } C_{bk}^{\text{out}} x_{bk} \label{eq:objFun_pde_BP} \\ 
    & \hspace{0.05cm} \text{ s.t. } \sum_{b \in \mathcal{B}_k} x_{bk} = y_k^{\text{p}},  \hspace{5.1cm}  k \in \mathcal{K}, \label{eq:periodic_bp}\\
    & \hspace{0.5cm} x_{bk} \leq y_k^{\text{p}} z_b, \hspace{5.6cm}  k \in \mathcal{K},  b \in \mathcal{B}_k,  \label{eq:selected_blocks_bp}\\
    & \hspace{0.5cm} v_{b}^{53} =  \max \left [ 0, \left \lceil \frac{1}{2} \left ( \sum_{k \in \mathcal{K}_b, \gamma_k = 53} x_{bk} - \sum_{k \in \mathcal{K}_b, \gamma_k=40 } x_{bk} \right ) \right \rceil \right ],    \nonumber \\
    & \hspace{7.6cm} b \in \mathcal{B}, \label{eq:cap1_bis}\\
    & \hspace{0.5cm} v_{b}^{40} = \left \lceil \frac{1}{2} \left ( \sum_{k \in \mathcal{K}_b} x_{bk}  \right ) \right \rceil - v_{b}^{53}, \hspace{2.95cm}  b \in \mathcal{B}, \label{eq:cap2_bis} \\
    & \hspace{0.5cm} \sum_{b \in \mathcal{B}_a} \left( L^{40} v_b^{40} + L^{53} v_b^{53} \right) \leq u_a, \hspace{3cm}  a \in \mathcal{A}^{TM},  \label{eq:cap3_bis} \\
    & \hspace{0.5cm} z_{b} \in \{0,1\}, \hspace{5.6cm}  b \in \mathcal{B},  \\
    & \hspace{0.5cm} v_b^{40}, v_b^{53} \in \mathbb{N}, \hspace{5.4cm}  b \in \mathcal{B}, \\
    & \hspace{0.5cm} x_{bk}  \in \mathbb{N}, \hspace{5.95cm}  k \in \mathcal{K}, b \in \mathcal{B}_k, \label{eq:pde_bp_end1}\\
    \textbf{wBP} & \hspace{0.6cm} \min_{x_1, \ldots, x_T} \sum_{t = 1}^T \sum_{k \in \mathcal{K}} \left[ \sum_{b \in \mathcal{B}_k \setminus \mathcal{B}_k^{\text{artif}} } C_{tbk}^{\text{flow}} x_{tbk} + \sum_{b \in \mathcal{B}_k^{\text{artif}} } C_{tbk}^{\text{out}} x_{tbk} \right] \label{eq:obj_pde_wbp}  \\ 
    & \hspace{0.6cm} \text{ s.t. } \sum_{b \in \mathcal{B}_k} x_{tbk} = y_{tk},  \hspace{4.3cm}  t \in \mathcal{T},  k \in \mathcal{K}, \label{eq:pde_wbp_2}\\
    & \hspace{1.2cm} x_{tbk} \leq y_{tk} z_b, \hspace{4.6cm}  t \in \mathcal{T},  k \in \mathcal{K},  b \in \mathcal{B}_k,  \label{eq:selected_blocks_wbp}\\
    & \hspace{1.2cm} v_{tb}^{53} =  \max \left [ 0, \left \lceil \frac{1}{2} \left ( \sum_{k \in \mathcal{K}_b, \gamma_k = 53} x_{tbk} - \sum_{k \in \mathcal{K}_b, \gamma_k=40 } x_{tbk} \right ) \right \rceil \right ], \nonumber\\
    & \hspace{7.55cm} t \in \mathcal{T},  b \in \mathcal{B}, \label{eq:bp_platf53_2} \\
    & \hspace{1.2cm} v_{tb}^{40} = \left \lceil \frac{1}{2} \left ( \sum_{k \in \mathcal{K}_b} x_{tbk}  \right ) \right \rceil - v_{tb}^{53}, \hspace{2.12cm}  t \in \mathcal{T},   b \in \mathcal{B},  \label{eq:bp_platf40_2}\\
    & \hspace{1.2cm} \sum_{b \in \mathcal{B}_a} \left( L^{40} v_{tb}^{40} + L^{53} v_{tb}^{53} \right) \leq u_a, \hspace{2.25cm}  t \in \mathcal{T},  a \in \mathcal{A}^{TM}, \label{eq:cap4_bis}  \\
    & \hspace{1.2cm} v_{tb}^{40}, v_{tb}^{53} \in \mathbb{N}, \hspace{4.58cm}  t \in \mathcal{T}, b \in \mathcal{B} \label{eq:last_1_P}, \\
    & \hspace{1.2cm} x_{tbk} \in \mathbb{N}, \hspace{5cm}  t \in \mathcal{T},  k \in \mathcal{K},  b \in \mathcal{B}_k. \label{eq:last_P}
\end{flalign}
\end{footnotesize}
Constraints~\eqref{eq:periodic_bp} and \eqref{eq:pde_wbp_2} enforce that the demand is satisfied by either the network capacity or outsourcing. Constraints~\eqref{eq:selected_blocks_bp} and \eqref{eq:selected_blocks_wbp} ensure flows to be on selected blocks only. Constraints~\eqref{eq:cap1_bis},~\eqref{eq:cap2_bis},~\eqref{eq:bp_platf40_2} and \eqref{eq:bp_platf53_2} fix the number of platforms required to transport the demand. These constraints model the double-stacking possibilities for the containers of different sizes. The 40-foot platforms are preferred  to 53-foot platforms, and the remaining 53-foot containers are stacked on top of 40-feet containers. The platform lengths are denoted $L^{40}$ and $L^{53}$, respectively. Constraints~\eqref{eq:cap3_bis} and \eqref{eq:cap4_bis} ensure that the train capacity is not exceeded. The latter, denoted $u_a, a\in \mathcal{A}^{TM}$, is defined for the set of arcs that represent moving trains, $\mathcal{A}^{TM}$. The set of blocks that use train moving arc $a \in \mathcal{A}^{TM}$ is $\mathcal{B}_a$.

The main costs from $C^{\text{PDE}}$ come from the flow of commodities on blocks and the outsourcing flows. The cost of building the blocks is relatively low compared to the containers being moved on said blocks. The value of the parameters in the objective function, namely $C^{\text{design}}, C^{\text{flow}}$ and $C^{\text{out}}$ was defined after a thorough analysis with CN of the paths produced by \textbf{BP} for each commodity.

\section{Results}
\label{section:results}

To assess the performance of the proposed methodology and algorithms, we solve instances that differ by the number of commodities, the resource sharing and the capacity constraints. We present below three metrics that allow to quantify the resource sharing and the capacity constraints. We report in Section~\ref{section:without_clustering} the results of the metaheuristics on the different instances without the first clustering heuristic step, and with the clustering heuristic in Section~\ref{section:with_clustering}.

Let us consider an illustration to explain the aforementioned metrics, with four commodities $k_1, k_2, k_3, k_4$ and two trains $A_1, A_2$. Commodities $k_1$ and $k_2$ take $A_1$ and $A_2$ on their route, $k_3$ takes $A_1$ and $k_4$ takes $A_2$. To quantify the resource sharing, we consider two metrics. We denote $\tau$ the average number of commodities that a train can carry. In the example, $\tau = 3$ since both trains carry a total of 3 commodities. We also introduce $\kappa$, a metric which quantifies, for each commodity, the number of other commodities sharing a train with the former. In other words, $\kappa$ is the average number of commodities with at least a common train, per commodity. In the example, $\kappa = 2.5$ since $k_1$ and $k_2$ share the capacity with 3 other commodities and $k_3$ and $k_4$ with 2.

To quantify the capacity limits, we analyze the potentially outsourced commodities. We first solve \textbf{BP} with the real life instance with $K=170$ commodities and $\boldsymbol{\alpha} = \mathbf{1}_K$, i.e, the periodic demand in \textbf{BP} is the average of the demand forecasts. 
We partition the set of commodities $\mathcal{K}$ in two subsets: the commodities with completely satisfied demand and the commodities with partial or total outsourced demand. We designate the latter by $\mathcal{K}_{L}$.  

Table~\ref{table:instance_demand} reports the characteristics of each instance. Each instance corresponds to a subset of CN's network and can in fact be considered as a smaller transport network on its own. They have different sizes, and present various levels of resource sharing and capacity constraints. We also indicate the number of blocks $B = |\mathcal{B}|$ in Table~\ref{table:instance_demand}, that is the number of paths on the space-time graph for each instance. The instance $IC$ corresponds to the real life instance. 
Instances $I1$ and $I3$ are not constrained by capacity, as $ \mathcal{K}_L= 0$, and instances $I4$ and $IC$ have large capacity sharing metrics.

\begin{table}[htbp!]
\centering
\begin{tabular}{llclccc}
Instance     & $|\mathcal{K}|$  & $|\mathcal{K}_{L}|$ & $B$ & $\tau$ & $\kappa$\\ \hline
$I1$  & 28  &  0   &  328  & 2 &   4 \\
$I2$  & 26  &  17  &  487  & 2 &   3\\
$I3$  & 48  &  0   &  501  & 3 &   4\\
$I4$ & 55  &  12  &  991  & 4 &   11 \\
$IC$ & 170 & 84 & 2208 & 9 & 22
\end{tabular}
\caption{Description of the instances and their characteristics}
\label{table:instance_demand}
\end{table}

\paragraph{Solving \textbf{BP} and \textbf{wBP}} In the remainder of the paper, we present results where \textbf{ePDE} is solved either with the metaheuristics NS and NSDI presented in Section~\ref{section:metaheuristics}, or with the off-the-shelf solver NOMAD. The three algorithms use the sequential property and \textbf{BP} and \textbf{wBP} are solved sequentially with the commercial solver CPLEX 12.10.0 for a fixed periodic demand. For our application, we need to solve them to near optimality (the gap stopping criteria is fixed to 0.4\%) to avoid large variations in the objective function. Such variations can happen because of changes in outsourcing. If the outsourcing is due to a poor quality solution, it can lead to an erroneous search direction.
Experiments were performed on a computer Intel(R) Core(TM) i9-10980XE CPU @ 3.00GHz with 8 CPUs and 128 GB de RAM.

\paragraph{Parameters for NS and NSDI} The metaheuristics rely on several parameters, namely the number of neighbors $V$, the distance parameter $\beta$, and for NSDI, the diversification and intensification parameters $b^-, b^+$, $v^+$ and the maximum number of iterations without improvements $M$. These parameters are specific to the basic characteristics of the transport network (e.g., terminals and OD pairs). In our case, each instance corresponds to a different network and hence we tuned them for each instance.
We tested various sets of parameters and report results obtained with best performing parameters. For all instances, $b^- = 0.7$, $b^+ = 1.3$, $v^+ = 1.1$. For instances $I1$ to $I4$, $V = 15$, $\beta = 0.05$ and $M = 15$. Finally, for instance $IC$, $V = 10$, $\beta = 0.02$ and $M = 7$.

\subsection{Results Without Clustering }
\label{section:without_clustering}

We define $\boldsymbol{\alpha}_{q3}$ such that $\mathbf{y}_{\text{q3}}^{\text{p}} = \boldsymbol{\alpha}_{\text{q3}} \odot \mathbf{y}_{\text{mean}}^{\text{p}}$, where $y_{\text{q3}, k}^{\text{p}} = Q_{0.75}( y_{1k},\ldots,y_{Tk} )$. The vector $\boldsymbol{\alpha}_{q3}$ is the vector of deviation from the third quartile of the demand forecasts to the mean of the forecasts.

For each instance, we solve the formulation \textbf{ePDE}~\eqref{eq:obj_pde}-\eqref{eq:last_P} with either $K$ variables, i.e. one $\alpha_k$ per commodity $k$, or one variable when the components of $\boldsymbol{\alpha}$ are constrained to be equal. We run the following experiments: 
\begin{itemize}
    \item \textbf{ePDE} solved with the vector $\boldsymbol{\alpha}$ constrained to be equal to $\boldsymbol{\alpha}_{\text{max}}$ and $\boldsymbol{\alpha}_{q3}$ respectively. The feasible space is in fact a set of cardinality 1. 
    \item \textbf{ePDE} solved with NS, NSDI and NOMAD. We note that for the instance $IC$, NOMAD cannot be used as the number of variables is too large. 
    \item $\textbf{ePDE}$ solved with all components of $\boldsymbol{\alpha}$ constrained to be equal to $\alpha_{\text{mean}} = 1$, \\$\alpha_{\text{max}} = \max \{ \alpha_{\text{max}, k}, k=1, \ldots, K\}$ or $\alpha_{q3} = \max \{ \alpha_{q3, k}, k=1, \ldots, K\}$. 
    \item $\textbf{ePDE}$ with all components of $\boldsymbol{\alpha}$ constrained to be equal solved with NS, NSDI and NOMAD.
\end{itemize}

Solving $\textbf{PDE}$ with $\boldsymbol{\alpha}$ fixed to either $\boldsymbol{\alpha}_{\text{max}}$, $\boldsymbol{\alpha}_{q3}$ or $\alpha_{\text{mean}}$ constitutes the enumeration approach proposed in \cite{Laage2021}. Table~\ref{table:small_instances_allresults} reports the gap to best known solution for the different experiments and for each instance. Several findings emerge. Defining the periodic demand as a deviation from the average of the demand values allows to obtain good solutions and can reduce substantially the costs compared to \cite{Laage2021} or the common practice of simply taking the mean of the demand values. 
Moreover, having a single variable $\alpha$ and hence a restricted feasible space of first-level variables allows to obtain good solutions, even the best found solution except for instance $I4$. 

\begin{table}[htbp!]
\centering
\begin{tabular}{ll|rrrrrr}
 &  & $I1$ & $I2$ & $I3$  & $I4$ & $IC$ \\
 \hline
\multirow{6}{*}{\begin{tabular}[c]{@{}l@{}}scalar $\alpha$  \\
\end{tabular}} & $\alpha_{\text{mean}}$ & 0\% & 125\% & 361\%  & 137\%  & 51\%  \\
 & $\alpha_{\text{max}}$ & 0\% & 80\% & 0\% & 1010\% & 885\% \\
 & $\alpha_{\text{q3}}$  & 0\% & 22\% & 8\% &  338\%  & 228\%  \\
 & NS & 0\%  & 0\% & 8\% &  14\%  & 16\%  \\
 & NSDI & 0\% & 0\% & 0\% &  3\%  & 0\%  \\
 & NOMAD & 0\% & 2\% & 0\% &  3\%  &  10\% \\ \hline
\multirow{5}{*}{\begin{tabular}[c]{@{}l@{}}$K$-vector  $\boldsymbol{\alpha}$ \\ 
\end{tabular}} & $\boldsymbol{\alpha_{\text{max}}}$ & 0\% & 0\% & 0\% &   210\%  &  78\% \\ 
 & $\boldsymbol{\alpha}_{\text{q3}}$ & 0\% & 42\% & 16\% &  27\%  &  21\% \\
 & NS & 0\% & 125\% & 119\% &  64\%  &  51\% \\
 & NSDI & 0\% & 0\% & 0\% &  0\%  & 51\%  \\
 & NOMAD & 0\% & 31\% & 16\% & 3\% & -  
\end{tabular}
\caption{Gap to best known value}
\label{table:small_instances_allresults}
\end{table}

When the capacity sharing metrics are low, for instances $I1$, $I2$ and $I3$, the best objective function value is found with the periodic demand  $\mathbf{y}^{\text{p}} = \boldsymbol{\alpha}_{\text{max}}$. In other words, when the network is not constrained, or slightly constrained by capacity, we should use the maximum of the demand values as periodic demand. For the instance $I1$, $C^{\text{PDE}}$ does not vary with the periodic demand. Having one block per commodity is sufficient, so as long as the periodic demand is not zero, \textbf{BP} creates one block per commodity. If the periodic demand is large, more blocks are built in \textbf{BP} but not used in \textbf{wBP}. The main part of the tactical costs comes from the flow costs, hence the increase of $C^{\text{PDE}}$ due to the higher number of blocks is not significant. The variations of $C^{\text{PDE}}$ are not visible in Table~\ref{table:small_instances_allresults} because of their small magnitude.

The differences in costs resulting from a single common $\alpha_{\text{max}}$ or a vector $\boldsymbol{\alpha}_{\text{max}}$ are due to the limited capacity. 
When \textbf{ePDE} is solved with $\alpha_{\text{max}}$, the periodic demand is overestimated for the commodities where $\alpha_{\text{max},k}$ < $\alpha_{\text{max}}$. As a result, \textbf{BP} builds too many blocks for one commodity, resulting in a lack of capacity for the others, and in turn, outsourcing. 

For the real-life size instance $IC$, the best performance is reached with a single common deviation coefficient $\alpha$. The metaheuristic NSDI finds this solution. With a vector of $K=170$ variables, both heuristics fail to even improve the initial solution.

\subsection{Results With Clustering}
\label{section:with_clustering}

In this section, we restrict the analysis to instances $I4$ and $IC$ as they both have high resource sharing and limited capacity metrics. In the other instances, the network has enough capacity and the best solution is obtained with $\boldsymbol{\alpha}_{\text{max}}$. 
We perform the clustering described in Section~\ref{section:clustering} and report in Table~\ref{tab:clusters} the number of clusters for each instance:
\begin{itemize}
    \item CV: Variance-based clustering
    \item CR: Service-based clustering
    \item CRU: Relaxed service-based clustering
\end{itemize}

\begin{table}[htbp!]
    \centering
    \begin{tabular}{lll}
            &   $I4$   & $IC$ \\ \hline
       CV   &   5      &  5   \\
       CR   &   4      &  12  \\
       CRU  &   5      &  16 
    \end{tabular}
    \caption{Number of clusters created in each clustering step}
    \label{tab:clusters}
\end{table}

For CV, we conducted several experiments varying the number of clusters, but report only the best results found with $n_C = 5$. The five clusters are following:
\begin{itemize}
    \item $C_1 = \{k \in \mathcal{K}, \sigma_{\text{min}}^{\text{s}} \leq \sigma_k^{\text{s}} \leq Q_{0.25} \} $
    \item $C_2 = \{k \in \mathcal{K}, Q_{0.25} < \sigma_k^{\text{s}} \leq  Q_{0.5} \} $
    \item $C_3 = \{k \in \mathcal{K}, Q_{0.5}  < \sigma_k^{\text{s}} \leq  Q_{0.75} \} $
    \item $C_4 = \{k \in \mathcal{K}, Q_{0.75} < \sigma_k^{\text{s}} \leq  Q_{0.9} \} $
    \item $C_5 = \{k \in \mathcal{K}, Q_{0.9}  < \sigma_k^{\text{s}} \leq  \sigma_{\text{max}}^{\text{s}} \} $, 
\end{itemize}
where $\sigma_{\text{min}}^{\text{s}} = \min \{\sigma_k^{\text{s}} , k=1,\ldots,K \}$, $\sigma_{\text{max}}^{\text{s}} = \max \{\sigma_k^{\text{s}} , k=1,\ldots,K \}$. Recall that $Q_e$ denotes the e-th percentile of $\{ \sigma_k^{\text{s}} , k=1,\ldots,K\}$.

Table~\ref{table:results_with_clusters} reports the gap to the best found solution after solving \textbf{ePDE} with the aforementioned clusters. For comparison, we include the results from Table~\ref{table:small_instances_allresults} in the first two sets of rows.
The results clearly show the benefit of the two-step heuristic as combining clustering and metaheuristics lead to reduced costs. In the case of $I4$ the costs are reduced by 2 percentage points, and by 5 percentage points for $IC$.

\begin{table}[H]
\centering
\begin{tabular}{ll|ll}
 &  & $I4$ & $IC$ \\ \hline
\multirow{6}{*}{\begin{tabular}[c]{@{}l@{}}scalar $\alpha$ \\ 
\end{tabular}} & $\alpha_{\text{mean}}$ & 143\%  & 59\%  \\
 & $\alpha_{\text{max}}$ & 1037\%  & 938\% \\
 & $\alpha_{\text{q3}}$ & 349\%   & 246\% \\
 & NS &  17\%  & 22\% \\
 & NSDI & 5\%  & 5\% \\
 & NOMAD & 5\%  & 16\% \\ \hline
\multirow{5}{*}{\begin{tabular}[c]{@{}l@{}}$K$-vector $\boldsymbol{\alpha}$ \\ 
\end{tabular}} & $\boldsymbol{\alpha_{\text{max}}}$ & 218\%  & 87\% \\ 
 & $\boldsymbol{\alpha_{\text{q3}}}$ & 30\%  & 27\% \\
 & NS & 68\%  & 59\% \\
 & NSDI & 2\%   & 59\% \\
 & NOMAD &  6\%  & - \\ \hline
 \multirow{3}{*}{\begin{tabular}[c]{@{}l@{}} CR 
 \end{tabular}} & NS & 1\%  & 59\% \\ 
 & NSDI & 1\%   & 2\% \\
 & NOMAD &  1\%  & \textbf{0\%} \\ \hline
 \multirow{3}{*}{\begin{tabular}[c]{@{}l@{}} CRU 
 \end{tabular}} & NS & 6\%  & 13\% \\ 
 & NSDI &  1\%  & 4\% \\
 & NOMAD &  2\%  & 2\% \\ \hline
  \multirow{3}{*}{\begin{tabular}[c]{@{}l@{}} CV 
  \end{tabular}} & NS & 21\% & 25\% \\ 
 & NSDI &  \textbf{0\%} & 12\% \\
 & NOMAD & \textbf{0\%}  & 1\% \\
\end{tabular}
\caption{Gap to best known value with clustering}
\label{table:results_with_clusters}
\end{table}

Clusters created from CV lead to a plan able to handle changes in demand at each period. They define a higher periodic demand for commodities whose demand forecasts have a large variance over the tactical planning periods. Then, when the forecasts in one period are higher than the average over all periods, the design is appropriate to handle the demand increase. 
Clusters from CR allow to increase the capacity for commodities competing for the same services, which lead to a reduction of the outsourcing. 

The best performance is obtained when solving \textbf{ePDE} 
with the clustering step for both the metaheuristics and NOMAD. Therefore, when clustering reduces the feasible space, it allows to leverage a state-of-the-art black-box optimization solver for large-scale applications and produces higher-quality solutions.

\begin{figure}
    \begin{subfigure}[b]{\linewidth}
        \includegraphics[width=\linewidth]{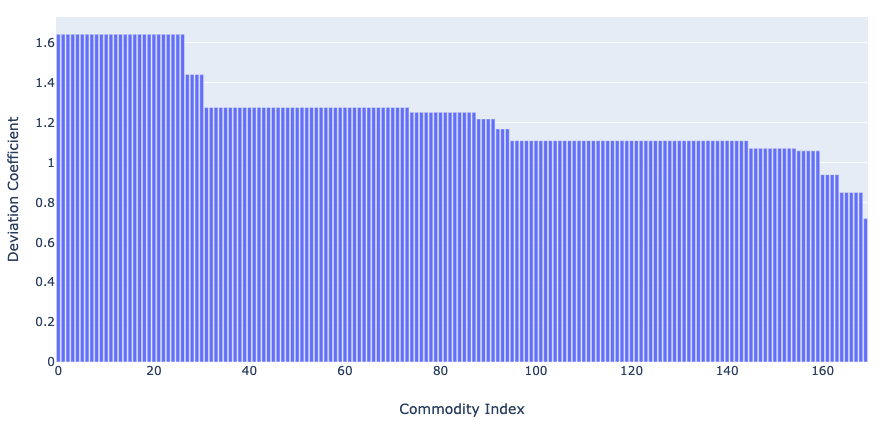}
        \caption{Solution with CR clustering and NSDI}
        \label{fig:nsdi_cr}  
    \end{subfigure}
    \begin{subfigure}[b]{\linewidth}
        \includegraphics[width=\linewidth]{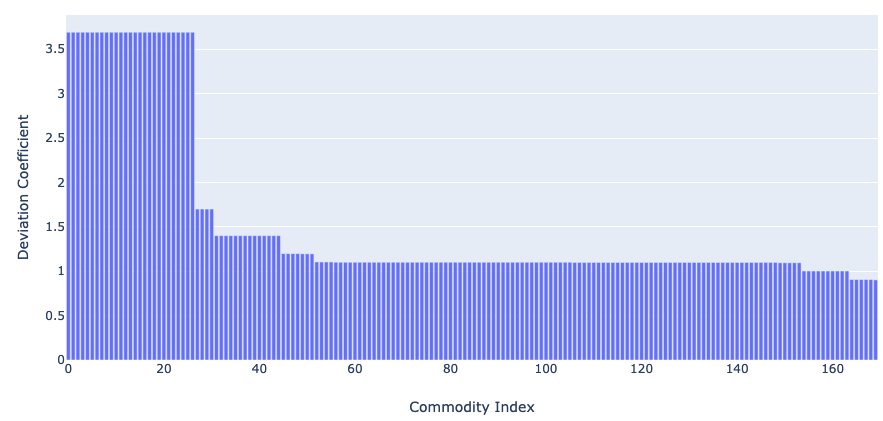}
        \caption{Solution with CR clustering and NOMAD}
        \label{fig:nomad_cr}
    \end{subfigure}
\caption{Values of the coefficient deviations of all commodities obtained after solving \textbf{ePDE}. Each bar corresponds to one of the 170 commodities in instance $IC$. }
\label{fig:example_alpha_solutions}
\end{figure}

To illustrate how the solutions can be interpreted, we present in Figure~\ref{fig:example_alpha_solutions} the values of the deviation coefficients over all commodities for CR combined with NSDI (a) and NOMAD (b). We can identify one cluster for which both algorithms define a periodic demand that is higher than the average demand forecast ($\alpha =1$). The solution from NOMAD for this cluster is larger ($\alpha = 3.69$) than the one from NSDI ($\alpha = 1.64$).
This is because NOMAD allows for more diversification in the first iterations. Recall that NOMAD finds a solution that is of better quality than NSDI.

We end this section with a discussion about computational cost. Since tactical planning concerns medium-term decisions, is not critical to solve \textbf{ePDE} in a short computing time. The computing budget depends on the application but would typically be in the range of a few hours (overnight computation) to a day. In our case, the time depends on the number of iterations of the metaheuristics and the time it takes to solve  \textbf{BP} and \textbf{wBP} at each iteration. The latter depends on the required optimality gap (and fixed time limit) but also exhibits a variation intrinsic to MIP solvers. When either \textbf{BP} or \textbf{wBP} successfully reached the optimality gap criteria at 0.4\%, they reached it on average in 1.97s and 5.06s respectively for instance $I4$ and in 109.1s and 85.8s for instance $IC$.
We used a 900s time limit or optimality gap of 0.4\% as stopping critera for either \textbf{BP} and \textbf{wBP}.  

We now focus our analysis on the number of evaluations required by the metaheuristics, that is the number of times the algorithms solve the lower levels \textbf{BP-wBP} before reaching the best found solution. We report the results in Table~\ref{table:computing_times}. 
The results show that while NOMAD finds the best solution, the method NSDI gives a good solution in less evaluations, hence less time. For the instance $IC$ with CR clustering, NSDI returns a solution that is 2 percentage points more expensive than NOMAD, but in less than a third of the evaluations. 

\begin{table}[htbp!]
    \centering
    \begin{tabular}{ll|ll}
 &  & $I4$ & $IC$ \\ \hline
\multirow{3}{*}{\begin{tabular}[c]{@{}l@{}}scalar $\alpha$ \\ 
\end{tabular}} 
 & NS & 31   &  46 \\
 & NSDI & 46 &  222 \\
 & NOMAD & 26 &  52 \\ \hline
\multirow{3}{*}{\begin{tabular}[c]{@{}l@{}}$K$-vector $\boldsymbol{\alpha}$ \\ 
\end{tabular}}
 & NS &  31  & 1 \\
 & NSDI &  79  & 1 \\
 & NOMAD &  168  & - \\ \hline
 \multirow{3}{*}{\begin{tabular}[c]{@{}l@{}} CR 
 \end{tabular}} 
 & NS & 31  & 1 \\ 
 & NSDI &  61  &  69 \\
 & NOMAD &   107  &  \textbf{225} \\ \hline
 \multirow{3}{*}{\begin{tabular}[c]{@{}l@{}} CRU 
 \end{tabular}} 
 & NS &  76  &  31 \\ 
 & NSDI &   107  &  247 \\
 & NOMAD &   126  &  450 \\ \hline
  \multirow{3}{*}{\begin{tabular}[c]{@{}l@{}} CV 
  \end{tabular}}
  & NS &  16 &  46 \\ 
 & NSDI &  \textbf{367}  &  11 \\
 & NOMAD &  \textbf{214}  &  184 \\
\end{tabular}
    \caption{Number of evaluations. We indicate by 1 when the initial solution is not improved, and by a dash when the algorithm could not solve the problem. The figures in bold correspond to the best known solution in Table~\ref{table:results_with_clusters}.}
    \label{table:computing_times}
\end{table}

\section{Conclusion}
\label{section:conclusion}

This paper builds on the work in \cite{Laage2021}. They introduced the Periodic Demand Estimation problem and showed that it is of high practical importance. However, for the sake of computational tractability, they made strong restrictive assumptions on the solution space. In this paper, we relaxed those assumptions and introduced a new formulation where the periodic demand is defined as a deviation from averaged time series forecasts. This results in a bounded continuous solution space where variables represent deviation coefficients. The latter have an intuitive interpretation as the deviation from average demand values. The new formulation remains challenging to solve because of the potentially large number of deviation coefficients (one per commodity).

We proposed a two-step heuristic. The objective of the first step is to reduce the number of variables in a data driven way through heuristic variable fixing. 
We proposed clustering approaches that group commodities for which the deviation coefficients are assumed equal. The new formulation and this first step constituted our main methodological contribution as they allow to leverage state-of-the-art metaheuristics in the second step (black-box optimization solver NOMAD or local search algorithms).

We reported results for a real large-scale application at the Canadian National Railway Company, that are improved compared to \cite{Laage2021}. 
By reducing the number of variables, the clustering step makes the use of NOMAD possible for this large-scale application which originally had 170 variables. 
While they did not find the best solution for all instances, the local search algorithms proposed in this paper allowed to reach solutions 1-2 percentage points more expensive than the best one, but at a substantially smaller computational cost than NOMAD.

Finally, we identify three main directions for future work. In this paper we limited the clustering to data-driven heuristics and the three algorithms had somewhat similar performance. In future work we plan to investigate unsupervised learning algorithms based on engineered features capturing demand and supply interactions. A second direction concerns the reduction of computational cost. The algorithm solving the PDE problem relies on accurate information about the lower-level problems. The latter are costly to solve. We plan to investigate ways to speed up the computation of the lower-level problems while maintaining the accuracy needed to guide the search for the PDE solution. This could be done by adequately relaxing the problems, reusing information about the lower-level problems from one iteration to the next, or predicting the solution values \citep[similar to][]{LarsenEtAl21}. The third direction for future research relates to stochastic programming. Our approach may be viewed as generating one single demand scenario. For the sake of computational tractability of large-scale problems, we focus on deterministic formulations. It would, however, be interesting to investigate if it is valuable and tractable to solve a sample average approximation of the problem based on a few, rather than just one scenario.

\section*{Acknowledgments}
We gratefully acknowledge the close collaboration with personnel from different divisions of the Canadian National Railway Company.
This research was funded by the Canadian National Railway Company Chair in Optimization of Railway Operations at Université de Montréal and a Collaborative Research and Development Grant from the Natural Sciences and Engineering Research Council of Canada (CRD-477938-14).

\end{document}